\documentclass{article}
\usepackage{amssymb,amsmath,amsthm}

\theoremstyle{plain}
\newtheorem{lemma}{Lemma}[section]
\newtheorem*{theorem}{Theorem}
\newtheorem{corollary}[lemma]{Corollary}

\newtheorem*{Zariski problem}{Zariski problem}
\newtheorem*{ZCP}{Zariski cancellation problem}

\theoremstyle{definition}

\newtheorem*{remarks}{Remarks}

\theoremstyle{remark}

\DeclareMathOperator{\trdeg}{trdeg}

\DeclareMathOperator{\ak}{AK}

\DeclareMathOperator{\Frac}{Frac}

\DeclareMathOperator{\ch}{char}

\DeclareMathOperator{\Spec}{Spec}

\DeclareMathOperator{\EXP}{EXP}

\newcommand{\K}{\mathbf{k}}
\newcommand{\f}{\varphi}

\providecommand{\bysame}{\makebox[3em]{\hrulefill}\thinspace}

\begin{document}
\title{On the rigidity of small domains}
\author{A. Crachiola\\
Department of Mathematics\\
Wayne State University\\
Detroit, MI 48202, USA\\
crach@math.wayne.edu\\
\\
L. Makar-Limanov\thanks{Supported by an NSA grant.}\\
Department of Mathematics\\
Wayne State University\\
Detroit, MI 48202, USA\\
Department of Mathematics \& Computer Science\\
Bar-Ilan University\\
Ramat-Gan 52900, Israel\\
lml@math.wayne.edu, lml@macs.biu.ac.il}

\date{January 21, 2004}
\maketitle
\begin{abstract}
Let $\K$ be an algebraically closed field of arbitrary
characteristic. Let $A$ be an affine domain over $\K$ with
transcendence degree 1 which is not isomorphic to $\K[x]$, and let
$B$ be a domain over $\K$. We show that the AK invariant
distributes over the tensor product of $A$ by $B$. As a
consequence, we obtain a generalization of the cancellation
theorem of S.~Abhyankar, P.~Eakin, and W.~Heinzer.


\end{abstract}
\section{Introduction}
All rings in this article are commutative with identity. For a
ring $A$, let $A^{[n]}$ denote the polynomial ring in $n$
indeterminates over $A$. Let $\K$ be a field of arbitrary
characteristic. Consider the well known

\begin{ZCP}
Let $V_1$ and $V_2$ be affine varieties over $\K$ such that $V_1
\times \K^n \cong V_2 \times \K^n$ for some positive integer $n$.
Does it follow that $V_1 \cong V_2$?
\end{ZCP}

Significant results on this problem were published in 1972.
S.~Abhyankar, P.~Eakin, and W.~Heinzer~\cite{AEH} answered the
question affirmatively for affine curves using algebraic methods
(see Corollary~\ref{C:AEH}). Also, M.~Hochster~\cite{hochster}
gave a negative answer by constructing a 4-dimensional
counterexample over the real numbers. Given this counterexample,
it is natural to seek some restriction on $V_1$ and $V_2$ under
which we may solve the problem. Because the example given by
Hochster requires the formally real property of the real numbers,
a natural restriction is that $\K$ be algebraically closed.
However, in 1989 W.~Danielewski~\cite{danielewski} provided a
2-dimensional counterexample over the complex numbers which lead
to a class of similar counterexamples~\cite{fieseler,wilkens}.
Another classical restriction on the Zariski cancellation problem
is that $V_2$ be affine space. An affirmative answer to this case
for surfaces was given by T.~Fujita, M.~Miyanishi, and
T.~Sugie~\cite{fujita,miyanishiS} for characteristic 0 fields, and
P.~Russell~\cite{russell} extended their result to fields of
arbitrary characteristic. Also, T.~Fujita and
S.~Iitaka~\cite{FujIit} solved the problem affirmatively for
varieties $V_i$ of any dimension over $\mathbf{C}$ in the case
when the logarithmic Kodaira dimension of $V_i$ is not $- \infty$.
Beyond this, the problem is still open.

In this article, we consider another perspective. The Zariski
cancellation problem can be posed algebraically as follows. Let
$A_1$ and $A_2$ be affine domains over $\K$. Does $A_1^{[n]} \cong
A_2^{[n]}$ imply $A_1 \cong A_2$? Viewing these polynomial rings
as tensor products $A_i \otimes_{\K} \K^{[n]}$, we can pose a more
general cancellation question. If $B$ is an algebra over $\K$ such
that $A_1 \otimes_{\K} B \cong A_2 \otimes_{\K} B$, under what
conditions can we conclude that $A_1 \cong A_2$? Of course we
still have the counterexamples due to Danielewski and Hochster.
However, in light of the positive result of Abhyankar, Eakin, and
Heinzer, in this article we shall study the ``small''
1-dimensional case of this more general cancellation question. For
us this means that $A_1$ and $A_2$ have transcendence degree 1
over $\K$.

A fruitful approach to understanding cancellation is to study
additive group actions on a variety.  Over characteristic 0
fields, this means studying locally nilpotent derivations on the
ring of regular functions. Over prime characteristic fields, we
can analogously consider locally finite iterative higher
derivations. One tool which has been found beneficial in the
characteristic 0 setting is the AK invariant, defined for a
variety as the subring of regular functions which remain invariant
under all additive group actions on the variety. The main goal of
this article is to prove the following theorem, which has
immediate consequences on the general question of cancellation,
including the theorem of Abhyankar, Eakin, and Heinzer.

\begin{theorem}
Let $\K$ be an algebraically closed field. Let $A$ be an affine
domain over $\K$ with $\trdeg_{\K}(A)=1$ which is not isomorphic
to $\K^{[1]}$. Let $B$ be a domain over $\K$. Then
\[
\ak(A \otimes_{\K} B) = \ak(A) \otimes_{\K} \ak(B).
\]
\end{theorem}

Aside from the important geometric motivation behind this result,
it is valuable from the ring theoretic perspective as a means for
studying the structure of certain rings. Because isomorphisms of
rings restrict to isomorphisms of their AK invariants, we can use
the AK invariant as a probe into the automorphism group of a ring.
For any ring $A$, to say that $\ak(A)=A$ is to say that there are
no (nontrivial) exponential maps on $A$. (This notion will be
explained in the sequel.) This in turn tells us that $A$ is
lacking a certain type of automorphism. In case $A$ is a
$\K$-algebra, it means there are no nontrivial actions of the
additive group $\K^+$ on $A$. We will call a ring $A$
\emph{rigid\/} if $\ak(A)=A$. In fact, all domains with
transcendence degree 1 over $\K$ are rigid, with the exception of
$\K^{[1]}$ (see Lemma~\ref{L:dim1}). So the main theorem of this
article is a statement on the rigidity of such domains. We see,
for instance, that a tensor product of two rigid transcendence
degree 1 domains will remain rigid. Also, if $A \otimes_{\K} B$ is
not rigid, where $A$ and $B$ are as in the theorem, the
exponential maps on this tensor product leave $A$ fixed like an
anchor around which elements of $B$ are moved.

The slogan for this article is therefore, ``small rigid domains
stay rigid''. It remains to study the question of rigidity for
larger domains.

\section{Exponential maps and the AK invariant}
Let us review some relevant notions with a view towards the
definition of the AK invariant.

Let $\K$ be a field of arbitrary characteristic and let $A$ be a
$\K$-algebra. Suppose $\f:A \to A^{[1]}$ is a $\K$-algebra
homomorphism. We write $\f = \f_t:A \to A[t]$ if we wish to
emphasize an indeterminate $t$. We say that $\f$ is an
\emph{exponential map on $A$\/} if it satisfies the following two
additional properties.
\renewcommand{\theenumi}{\roman{enumi}}
\begin{enumerate}
\item $\varepsilon_0 \f_t$ is the identity on $A$, where
$\varepsilon_0:A[t] \to A$ is evaluation at $t=0$.\label{exp1}

\item $\f_s \f_t = \f_{s+t}$, where $\f_s$ is extended by
$\f_s(t)=t$ to a homomorphism $A[t] \to A[s,t]$.\label{exp2}
\end{enumerate}
(When $A$ is the coordinate ring of an affine variety $\Spec(A)$
over $\K$, the exponential maps on $A$ correspond to algebraic
actions of the additive group $\K^+$ on $\Spec(A)$~\cite[\S
9.5]{essen}.)

Given an exponential map $\f: A \to A[t]$, set $\f(t) = t$ to
obtain an automorphism of $A[t]$ with inverse $\f_{-t}$. Consider
the map $\varepsilon_1 \f :A \to A$, where $\varepsilon_1: A[t]
\to A$ is evaluation at $t=1$. One can check that $\varepsilon_1
\f$ is an automorphism of $A$ with inverse $\varepsilon_1
\f_{-t}$.

Define
\[
A^{\f} = \{ a \in A \,|\, \f(a)=a \},
\]
a subalgebra of $A$ called the \emph{ring of $\f$-invariants}.

For each $a \in A$ and each natural number $i$, let $D^i(a)$
denote the $t^i$-coefficient of $\f(a)$. Let $D = \{ D^0, D^1,
D^2, \ldots \}$. To say that $\f$ is a $\K$-algebra homomorphism
is equivalent to saying that the sequence $\{D^i(a)\}$ has
finitely many nonzero elements for each $a \in A$, that $D^n:A \to
A$ is $\K$-linear for each natural number $n$, and that the
Leibniz rule
\[
D^n(ab) = \sum_{i+j=n} D^i(a) D^j(b)
\]
holds for all natural numbers $n$ and all $a,b \in A$. The above
properties (\ref{exp1}) and (\ref{exp2}) of the exponential map
$\f$ translate into the following properties of $D$.
\renewcommand{\theenumi}{\roman{enumi}}
\begin{enumerate}
\item $D^0$ is the identity map.

\item (iterative property) For all natural numbers $i,j$,
\[
D^i D^j = \binom{i+j}{i} D^{i+j}.
\]
\end{enumerate}

Due to all of these properties, the collection $D$ is called a
\emph{locally finite iterative higher derivation on $A$}. More
generally, a \emph{higher derivation on $A$\/} is a collection $D
= \{D^i\}$ of $\K$-linear maps on $A$ such that $D^0$ is the
identity and the above Leibniz rule holds. The notion of higher
derivations is due to H.~Hasse and F.K.~Schmidt~\cite{hasse}.
Every higher derivation on $A$ has a unique extension to a higher
derivation on any given localization of $A$, determined through
extension of the Leibniz rule to fractions.~\cite[\S
27]{matsumura} When the characteristic of $A$ is 0, each $D^i$ is
determined by $D^1$, which is a locally nilpotent derivation on
$A$. In this case, $\f = \exp(tD^1) = \sum_i \frac{1}{i!}(tD^1)^i$
and $A^{\f}$ is the kernel of $D^1$.

Let $\EXP(A)$ denote the set of all exponential maps on $A$. We
define the \emph{AK invariant}, or \emph{ring of absolute
constants of $A$\/} as
\[
\ak(A) = \bigcap_{\f \in \EXP(A)} A^{\f}.
\]
This is a subalgebra of $A$ which is isomorphism preserved.
Indeed, any isomorphism $f:A \to B$ of $\K$-algebras restricts to
an isomorphism $f: \ak(A) \to \ak(B)$. To understand this, observe
that if $\f \in \EXP(A)$ then $f \f f^{-1} \in \EXP(B)$. We say
that $A$ is \emph{rigid\/} if $\ak(A) = A$. That is, the only
exponential map on $A$ is the standard inclusion $\f(a) = a$ for
all $a \in A$.

The above discussion of exponential maps, locally finite iterative
higher derivations, and the AK invariant makes sense more
generally for any (not necessarily commutative) ring. However, we
will not need this generality.

Given an exponential map $\f$ on a domain $A$ over $\K$, we can
define the \emph{$\f$-degree\/} of an element $a \in A$ by
$\deg_{\f}(a) = \deg_t(\f(a))$ (where $\deg_t(0) = - \infty$).
Note that $A^{\f}$ consists of all elements of $A$ with
non-positive $\f$-degree. The function $\deg_{\f}$ is a degree
function on $A$, i.e. it satisfies these two properties for all
$a,b \in A$.
\renewcommand{\theenumi}{\roman{enumi}}
\begin{enumerate}
\item $\deg_{\f}(ab) = \deg_{\f}(a) + \deg_{\f}(b)$.

\item $\deg_{\f}(a+b) \leq \max \{\deg_{\f}(a),\deg_{\f}(b) \}$.
\end{enumerate}

Equipped with these notions, we now collect some useful facts.

\begin{lemma}\label{L:facts}
Let $\f$ be an exponential map on a domain $A$ over $\K$. Let $D =
\{D^i\}$ be the locally finite iterative higher derivation
associated to $\f$.

\renewcommand{\theenumi}{\alph{enumi}}
\begin{enumerate}
\item If $a,b \in A$ such that $ab \in A^{\f} \setminus 0$, then
$a,b \in A^{\f}$.\label{factsa}

\item $A^{\f}$ is algebraically closed in $A$.\label{factsb}

\item For each $a \in A$, $\deg_{\f}(D^i(a)) \leq \deg_{\f}(a) -
i$. In particular, if $a \in A \setminus 0$ and $n =
\deg_{\f}(a)$, then $D^n(a) \in A^{\f}$.\label{factsc}
\end{enumerate}
\end{lemma}

\begin{proof}
(\ref{factsa}): We have $0 = \deg_{\f}(ab) = \deg_{\f}(a) +
\deg_{\f}(b)$, which implies that
$\deg_{\f}(a) = \deg_{\f}(b) = 0$.\\
(\ref{factsb}): If $a \in A \setminus 0$ and $c_n a^n + \cdots +
c_1 a + c_0 = 0$ is a polynomial relation with minimal possible
degree $n \geq 1$, where each $c_i \in A^{\f}$ with $c_0 \ne 0$,
then $a(c_n a^{n-1} + \cdots +c_1) = -c_0 \in A^{\f}
\setminus 0$. By part (\ref{factsa}), $a \in A^{\f}$.\\
(\ref{factsc}): Use the iterative property of $D$ to check that
$D^j(D^i(a))=0$ whenever $j > \deg_{\f}(a) - i$.
\end{proof}

\begin{lemma}\label{L:facts2}
Let $\f$ be a nontrivial exponential map (i.e not the standard
inclusion) on a domain $A$ over $\K$ with $\ch(A)=p \geq 0$. Let
$x \in A$ with minimal positive $\f$-degree $n$.

\renewcommand{\theenumi}{\alph{enumi}}
\begin{enumerate}
\item $D^i(x) \in A^{\f}$ for each $i \geq 1$. Moreover,
$D^i(x)=0$ whenever $i \geq 1$ is not a power of
$p$.\label{facts2a}

\item If $a \in A \setminus 0$, then $n$ divides
$\deg_{\f}(a)$.\label{facts2b}

\item Let $c = D^n(x)$. Then $A$ is a subalgebra of
$A^{\f}[c^{-1}][x]$, where $A^{\f}[c^{-1}] \subseteq
\Frac(A^{\f})$ is the localization of $A^{\f}$ at
$c$.\label{facts2c}
\end{enumerate}
\end{lemma}

\begin{proof}
In proving parts (\ref{facts2a}) and (\ref{facts2b}) we will
utilize the following fact. If $p$ is prime and $i = p^j q$ for
some natural numbers $i,j,q$, then $\binom{i}{p^j} \equiv q
\pmod{p}$~\cite[Lemma
5.1]{isaacs}.\\
(\ref{facts2a}): By part (\ref{factsc}) of Lemma~\ref{L:facts},
$D^i(x) \in A^{\f}$ for all $i \geq 1$. If $p=0$ then $n=1$, for
given any element in $A \setminus A^{\f}$ we can find an element
with $\f$-degree 1 by applying the locally nilpotent derivation
$D^1$ sufficiently many times. In this case, the second statement
is immediate. Suppose now that $p$ is prime and that $i>1$ is not
a power of $p$, say $i = p^j q$, where $j$ is a nonnegative
integer and $q \geq 2$ is an integer not divisible by $p$. Then
$D^{i - p^j}(x) \in A^{\f}$ and
\[
0 = D^{p^j} D^{i-p^j}(x) = \binom{i}{p^j}D^i(x) = q D^i(x).
\]
We can divide by $q$ to conclude that $D^i(x)=0$.\\
(\ref{facts2b}): Again if $p=0$ then $n=1$ and the claim is
obvious. Assume that $p$ is prime. By part (\ref{facts2a}) we have
$n = p^m$ for some integer $m \geq 0$. If $m=0$, the claim is
immediate. Assume that $m>0$. Let $d = \deg_{\f}(a)$. Suppose that
$p$ does not divide $d$. By part (\ref{factsc}) of
Lemma~\ref{L:facts}, $\deg_{\f}(D^{d-1}(a)) \leq 1$. Now,
$D^1D^{d-1}(a) = d D^d(a) \ne 0$. So $\deg_{\f}(D^{d-1}(a)) = 1 <
n$, contradicting the minimality of $n$. Hence we can write $d =
p^k d_1$ with $k \geq 1$ and $d_1$ not divisible by $p$. Making a
similar computation, $D^{p^k} D^{d-p^k}(a) = d_1 D^d(a) \ne 0$.
This implies that $\deg_{\f}(D^{d-p^k}(a)) = p^k$. Since $n = p^m$
is minimal,
we must have $k \geq m$, and so $n$ divides $d$.\\
(\ref{facts2c}): Let $a \in A \setminus 0$. By part
(\ref{facts2b}) we can write $\deg_{\f}(a) = ln$ for some natural
number $l$. If $l = 0$ then $a \in A^{\f}$ and we are done. We use
induction on $l>0$. Elements $c^l a$ and $D^{ln}(a) x^l$ both have
$\f$-degree $ln$. Let us check that $D^{ln}(c^l a) =
D^{ln}(D^{ln}(a) x^l)$. First, $D^{ln}(c^l a) = c^l D^{ln}(a)$ by
the Leibniz rule and because $c^l$ is $\f$-invariant. Secondly,
since $D^{ln}(x^l) = D^n(x)^l = c^l$ and $D^{ln}(a)$ is
$\f$-invariant, we see that $D^{ln}(D^{ln}(a) x^l) = c^l
D^{ln}(a)$ as well. (Remark: Though the equality $D^{ln}(x^l) =
D^{n}(x)^l$ does follow from the Leibniz rule, it may be more
immediately observed as follows. $D^{n}(x)$ is the leading
$t$-coefficient of $\f(x)$, and $\f$ is a homomorphism. Hence the
leading $t$-coefficient of $\f(x^l)$ is also that of $\f(x)^l$.)
Therefore, the element $y = c^l a - D^{ln}(a)x^l$ has $\f$-degree
less than $ln$ and hence less than or equal to $(l-1)n$. By the
inductive hypothesis, $y \in A^{\f}[c^{-1}][x]$. So $a = c^{-l}(y
+ D^{ln}(a) x^l) \in A^{\f}[c^{-1}][x]$.
\end{proof}

\begin{lemma}\label{L:dim1}
Let $A$ be a domain over $\K$ with $\trdeg_{\K}(A)=1$. Then
$\ak(A) = \K$ if and only if $A \cong \K^{[1]}$. Otherwise,
$\ak(A)=A$.
\end{lemma}

\begin{proof}
To see that $\ak(\K[X]) = \K$ (where $X$ is an indeterminate),
observe that $\psi(X) = X + t$ defines an exponential map on
$\K[X]$ with ring of $\psi$-invariants $\K$. Suppose $\ak(A) \ne
A$. Let $\f \in \EXP(A)$ be nontrivial. Part (\ref{factsb}) of
Lemma~\ref{L:facts} implies that $A^{\f} = \K$. By part
(\ref{facts2c}) of Lemma~\ref{L:facts2}, $A \subseteq \K[x]$ for
some $x \in A$ with minimal positive $\f$-degree. So $A = \K[x]$.
\end{proof}

Thus $\K^{[1]}$ is the only trancendence degree 1 domain which is
not rigid. By considering exponential maps of the form $\f_i(X_j)
= X_j + \delta_{ij}t$, where $\delta_{ij}$ is the Kronecker delta,
one can see that $\ak(\K^{[n]}) = \K$ for each natural number $n$.
However, if $A$ is a domain with transcendence degree $n \geq 2$
over $\K$, then $\ak(A) = \K$ does not imply that $A \cong
\K^{[n]}$~\cite{bandman}. One example will be given in the next
section.

\section{The main result and corollaries}

For most of our statements, tensor products are of $\K$-algebras
over $\K$, and transcendence degrees are taken over $\K$. So we
write $\otimes$ and $\trdeg$ rather than $\otimes_{\K}$ and
$\trdeg_{\K}$. If we need to specify a different field, we will
decorate the notation.

We can extend any $\f \in \EXP(A)$ to an exponential map on $A
\otimes B$ by defining $\f(b) = b$ for all $b \in B$. In other
words, we set $\f(\sum_i a_i \otimes b_i) = \sum_i \f(a_i) \otimes
b_i$. Any element of $\ak(A \otimes B)$ must be invariant under
such exponential maps, and so $\ak(A \otimes B) \subseteq \ak(A)
\otimes B$. We can interchange the roles of $A$ and $B$ in that
argument and further conclude that $\ak(A \otimes B) \subseteq
\ak(A) \otimes \ak(B)$.

\begin{theorem}\label{T:main}
Let $\K$ be an algebraically closed field. Let $A$ be an affine
domain over $\K$ with $\trdeg(A)=1$ which is not isomorphic to
$\K^{[1]}$. Let $B$ be a domain over $\K$. Then
\[
\ak(A \otimes B) = \ak(A) \otimes \ak(B).
\]
\end{theorem}

Of course, this can also be written as $\ak(A \otimes B) = A
\otimes \ak(B)$. The conclusion is false when $A \cong \K^{[1]}$,
as discussed below. This theorem has some immediate corollaries.

\begin{corollary}\label{C:1}
Let $\K$ be an algebraically closed field. Let $A_1$ and $A_2$ be
affine domains over $\K$ with $\trdeg(A_i) = 1$, $i=1,2$. Let $B$
be a domain over $\K$ such that $\ak(B)=\K$. If $A_1 \otimes B
\cong A_2 \otimes B$, then $A_1 \cong A_2$.
\end{corollary}

\begin{proof}
Suppose neither $A_1$ nor $A_2$ is isomorphic to $\K^{[1]}$.
Applying $\ak$ to both sides gives us $A_1 \otimes \K \cong A_2
\otimes \K$, so that $A_1 \cong A_2$. If $A_1 \cong \K^{[1]}$ but
$A_2 \ncong \K^{[1]}$, then
\[
A_2 \cong \ak(A_2 \otimes B) \cong \ak(A_1 \otimes B) \subseteq
\ak(A_1) \otimes \ak(B) \cong \K.
\]
But this is absurd.
\end{proof}

As a special case, take $B = \K^{[n]}$. The theorem implies that
$\ak(A^{[n]}) = \ak(A)$ for any affine domain with transcendence
degree 1 over $\K$. Since any isomorphism $f: A \to B$ restricts
to an isomorphism $f: \ak(A) \to \ak(B)$, we recover the
cancellation theorem of S.~Abhyankar, P.~Eakin, and W.~Heinzer:

\begin{corollary}[see \cite{AEH}]\label{C:AEH}
Let $\K$ be an algebraically closed field. Let $A_1$ and $A_2$ be
affine domains over $\K$ with $\trdeg(A_i) = 1$, $i=1,2$. If
$A_1^{[n]} \cong A_2^{[n]}$, then $A_1 \cong A_2$. Moreover, if
$f: A_1^{[n]} \to A_2^{[n]}$ and $A_1 \ncong \K^{[1]}$, then $f$
restricts to an isomorphism of $A_1$ onto $A_2$.
\end{corollary}

We feel compelled to again mention the geometric content of
Corollary~\ref{C:AEH}.

\begin{corollary}
Let $\K$ be an algebraically closed field. Let $V_1$ and $V_2$ be
affine curves over $\K$. If $V_1 \times \K^n \cong V_2 \times
\K^n$, then $V_1 \cong V_2$.
\end{corollary}

\begin{proof}
In algebraic terms it means $A_1^{[n]} \cong A_2^{[n]}$, where
$V_i = \Spec(A_i)$, $i=1,2$, and we must check that $A_1 \cong
A_2$.
\end{proof}

\begin{remarks} (1) The conclusion of Corollary~\ref{C:1} is still true for some
more general choices of $B$. For example, if $\ak(B) \cong
\K^{[n]}$ and $A_1 \otimes B \cong A_2 \otimes B$, then we can
apply $\ak$ twice to find that $A_1 \cong A_2$. (There are several
surfaces known to have AK invariant $\K^{[1]}$ or
$\K^{[2]}$~\cite{kaliman}.) For any choice of $B$ with finite
transcendence degree, we can apply $\ak$ to $A_i \otimes B$
several times in an attempt to show cancellation. Through each
application of $\ak$, the transcendence degree of the second
factor of the tensor product will decrease, unless it is rigid. So
$B$ must be rigid in any cancellation
counterexample with minimal dimension.\\
(2) As mentioned earlier, Corollary~\ref{C:AEH} is false when we
increase the transcendence degree of $A_1$ and $A_2$. The
following example is due to W.~Danielewski \cite{danielewski}. Let
$A_n$ be the coordinate ring of the surface $x^n y = z^2 -1$ over
the complex numbers $\mathbf{C}$. Then $A_1 \ncong A_2$ while
$A_1^{[1]} \cong A_2^{[1]}$. In fact, $A_i \ncong A_j$ whenever $i
\ne j$, but $A_i^{[1]} \cong A_j^{[1]}$ for all
$i,j$~\cite{fieseler,wilkens}. These domains also provide a
counterexample to the formula of our main result when $A \cong
\mathbf{C}^{[1]}$. One can prove that $\ak(A_1) = \mathbf{C}$ and
$\ak(A_2) = \mathbf{C}[x]$~\cite{ML2,ML3}. Now, $\ak(A_2^{[1]})
\cong \ak(A_1^{[1]}) \subseteq \ak(A_1) = \mathbf{C}$. Thus
$\ak(\mathbf{C}^{[1]} \otimes A_2) = \mathbf{C}$ while
$\ak(\mathbf{C}^{[1]}) \otimes \ak(A_2) = \mathbf{C}[x]$.\\
(3) In the special case $B = \K^{[n]}$, we can extend the theorem
(and hence its corollaries) to some non-algebraically closed
fields. Suppose $\mathbf{F}$ is a perfect field. (In particular,
$\mathbf{F}$ could be any characteristic 0 field.) Let $\K$ be an
algebraic closure of $\mathbf{F}$. Let $A$ be an affine domain
over $\mathbf{F}$ with $\trdeg_{\mathbf{F}}(A)=1$. It is known
that if $A \otimes_{\mathbf{F}} \K \cong \K^{[1]}$ then $A \cong
\mathbf{F}^{[1]}$~\cite{asanuma}. Using this fact we can easily
check that $\ak(A^{[n]}) = \ak(A)$ by considering the extension of
scalars $A \otimes_{\mathbf{F}} \K$ and applying the theorem.
\end{remarks}
\section{Proof of the main result}
Let $\K$ be an algebraically closed field. Let $A$ be an affine
domain over $\K$ with $\trdeg(A)=1$ which is not isomorphic to
$\K^{[1]}$. By Lemma~\ref{L:dim1}, $\ak(A) = A$. Let $B$ be a
domain over $\K$. We will view $A$ and $B$ as subalgebras of $A
\otimes B$ in the natural way. It is well known that the tensor
product of two affine domains is again an affine
domain~\cite{hartshorne}. Now suppose $z \in A \otimes B$ is a
zero divisor. Write $z = \sum_{i=1}^m a_i \otimes b_i$. Then $z$
belongs to the affine domain $A \otimes \K[b_1,\ldots,b_m]$. By
Theorem 36, Chapter III, of Zariski and Samuel~\cite{zariski}, $z$
is a zero divisor of this subdomain, and so $z=0$. Therefore, $A
\otimes B$ is a domain, and the lemmas of the previous section
apply to it.

Let us note that this next lemma does not require that $A$ have
transcendence degree 1. It still holds true for an affine domain
$A$ of any (necessarily finite) transcendence degree.

\begin{lemma}
If $A \subseteq \ak(A \otimes B)$ then $\ak(A \otimes B) = \ak(A)
\otimes \ak(B)$.
\end{lemma}

\begin{proof}
We need to show that $\ak(B) \subseteq \ak(A \otimes B)$. Let $b
\in \ak(B)$ and suppose that $\f(b) \ne b$ for some $\f \in \EXP(A
\otimes B)$. Let $f \in A \otimes B$ denote the leading
$t$-coefficient of $\f(b)$. Write $f = \sum_m a_m \otimes b_m$,
where the set $\{b_m\}$ is linearly independent over $\K$. Let
$\{g_1, \ldots, g_n\}$ be a finite generating set of $A$ over
$\K$. Since $\K$ is infinite, there exists a choice of values $c_i
\in \K$ for each $g_i$ such that evaluation of $g_i$ at $c_i$ is a
well-defined homomorphism whose kernel does not include the
element $a_1$. (In other words, there exists a point
$(c_1,\ldots,c_n) \in \Spec(A)$ which is not a zero of the regular
function $a_1$.) Let $\sigma: A \otimes B \to B$ denote the map
which sends $g_i \in A$ to $c_i$, $i=1,\ldots,n$, leaving all
elements of $B$ fixed. Let $\psi = \sigma \f$. We claim that $\psi
\in \EXP(B)$. It is clear that $\psi$ is a $\K$-homomorphism and
that $\varepsilon_0 \psi$ is the identity on $B$. Thus the
$t^i$-coefficients define a locally finite higher derivation
$\{\sigma D^i\}$ on $B$, and it remains to check that this
derivation is iterative. This follows routinely from the iterative
property of the higher derivation associated to $\f$ along with
the fact that $\f(a)=a$ for all $a \in A$. We leave the details to
the reader. Now $\psi$ is an exponential map on $B$, but the
$\psi$-degree of $b$ is larger than 0, since we chose $\sigma$ so
that $\sigma(f) \ne 0$. This contradicts our assumption that $b
\in \ak(B)$.
\end{proof}

So to prove the theorem we must demonstrate that $A \subseteq
\ak(A \otimes B)$. Suppose it is not the case. The next several
lemmas will bring this to a contradiction.

\begin{lemma}\label{L:technical1}
$A$ can be embedded in $\K^{[1]}$.
\end{lemma}

\begin{proof}
Since $A \nsubseteq \ak(A \otimes B)$, there exists $\f \in \EXP(A
\otimes B)$ for which $A \nsubseteq (A \otimes B)^{\f}$. In fact
then $A \cap (A \otimes B)^{\f} = \K$ by part (\ref{factsb}) of
Lemma~\ref{L:facts}. By part (\ref{facts2c}) of
Lemma~\ref{L:facts2}, $A$ is a subalgebra of $\Frac((A \otimes
B)^{\f})[x]$ for some $x \in A \otimes B$. Let $\{g_i\}$ be a
finite set of generators for $A$ over $\K$. Each $g_i$ is a
polynomial in $x$ with coefficients in $\Frac((A \otimes
B)^{\f})$. Let $\mathbf{E}$ be the subfield of $\Frac((A \otimes
B)^{\f})$ generated by all of these coefficients. Then $A$ is a
subalgebra of $\mathbf{E}[x]$ which is not contained in
$\mathbf{E}$. Since $\mathbf{E}$ is finitely generated over $\K$,
we can write $\mathbf{E} = \K(t_1,\ldots,t_k)[\alpha_1, \ldots,
\alpha_l]$ where elements $t_1,\ldots, t_k$ are transcendental
over $\K$, $\alpha_1$ is algebraic over $\K(t_1,\ldots,t_k)$, and
$\alpha_i$ is algebraic over $\K(t_1,\ldots,t_k)[\alpha_1, \ldots,
\alpha_{i-1}]$ for each $i = 2, \ldots, l$. We will choose values
in $\K$ for $t_1,\ldots,t_k$ and $\alpha_1, \ldots, \alpha_l$ in
order to embed $A$ in $\K[x]$. For this specialization process to
work successfully, we must choose these values so that certain bad
situations do not occur. First, we must insure that denominators
of elements of $A$ are not sent to zero. Next, we must insure that
the relations satisfied by each $\alpha_i$ do not become
contradictory (e.g. $0=1$). Finally, we do not want all elements
of $A$ to become elements of $\K$. Let us now outline how to avoid
these bad situations.

Each $g_i$ is a polynomial in $x$ with coefficients in
$\mathbf{E}$. Consider these coefficients as fractions with
numerators in $\K[t_1,\ldots,t_k][\alpha_1, \ldots, \alpha_l]$ and
denominators in $\K[t_1,\ldots,t_k]$. Let $d$ be a common
denominator for all of these fractions.  We will require $d$ to
not become 0.

Each $\alpha_i$ satisfies a polynomial $P_i$ with coefficients in
$\K[t_1,\ldots,t_k]$. Let $\overline{P_i}$ denote the polynomial
with coefficients in $\K$ that will be obtained by specializing
the coefficients of $P_i$. We will choose the values of
$t_1,\ldots,t_k$ so that the coefficients of $\overline{P_i}$ are
nonzero. Then we can chose a specialization value for $\alpha_i$
to be a root of $\overline{P_i}$ (since $\K$ is algebraically
closed). In this way, the relations satisfied by each $\alpha_i$
will not become contradictions. Let $p_i$ denote the product of
the coefficients of $P_i$ ($i=1,\ldots,l$). As long as $p_i$ does
not become zero, none of the coefficients of $P_i$ can become
zero.

We can assume that the generator $g_1 \in \mathbf{E}[x]$ has
$x$-degree at least 1. We want $g_1$ to specialize to a
nonconstant polynomial in $\K[x]$. The leading $x$-coefficient of
$g_1$ satisfies a polynomial $P_0$ with coefficients in
$\K[t_1,\ldots,t_k]$. Let $p_0$ be the constant coefficient of
$P_0$. We can assume $p_0 \ne 0$. If we specialize in a way that
$p_0$ does not become 0, then the leading $x$-coefficient of $g_1$
will not become 0, and hence $g_1$ will specialize to a
nonconstant polynomial.

Since $\K$ is infinite, there exists a $k$-tuple
$(c_1,\ldots,c_k)$ which is not a zero of the product $d p_0
\cdots p_l \in \K[t_1,\ldots,t_k]$, and hence not a zero of any of
its factors. Choose the value $c_i$ for $t_i$, $i=1,\ldots,k$. By
virtue of the above discussion, this choice results in a
$\K$-homomorphism $A \to \K[x]$ which sends $g_1$ to a nonconstant
polynomial. Since the image of this map has transcendence degree 1
over $\K$, the kernel must be trivial.
\end{proof}

Viewing $A$ as a subalgebra of $\K[x]$, L\"{u}roth's
theorem~\cite{isaacs} implies that the fraction field $\Frac(A)$
of $A$ is a simple extension of $\K$. Combining this with the
above lemma, we can choose the generator of $\Frac(A)$ over $\K$,
say $y$, to be a polynomial in $x$. From this it follows that $A
\subseteq \K[y]$. We will carry this assumption through the
remainder of the proof.

\begin{lemma}\label{L:technical2}
There exists a nonzero element $u \in A$ with the following
properties:
\begin{enumerate}
\renewcommand{\theenumi}{\roman{enumi}}
\item The ideal $\K[y]u$ of $\K[y]$ is contained in $A$.\label{u1}

\item If $a \in A$ such that $\K[y]a \subseteq A$, then $u$
divides $a$ in $\K[y]$.\label{u2}

\item If $f \in A \otimes B$ such that the ideal $(\K[y] \otimes
B)f$ of $\K[y] \otimes B$ is contained in $A \otimes B$, then $u$
divides $f$ in $\K[y] \otimes B$.\label{u3}
\end{enumerate}
\end{lemma}

\begin{proof}
Write $y = g h^{-1}$ where $g, h \in A$ and $h$ is a monic
polynomial in $y$ with degree $n$. Let us check that
$\K[y]h^{n-1}$ is contained in $A$. For $m=0,\ldots,n-1$ we have
$y^m h^{n-1}= g^m h^{n-m-1} \in A$. Now let $m \geq n$ and suppose
that $y^l h^{n-1} \in A$ for $0 \leq l <m$. Since $h = y^n +$
(terms of degree less than $n$), we can write $y^m = y^{m-n}h +
p_m(y)$, where $p_m(y)$ has degree at most $m-1$. Then $y^m
h^{n-1} = y^{m-n}h^n + p_m(y)h^{n-1} \in A$ by assumption. By
induction, $y^m h^{n-1} \in A$ for all $m \geq 0$. So the ideal
$\K[y]h^{n-1}$ is contained in $A$. Let $\mathfrak{a}$ be the
(nontrivial) ideal generated by all ideals of $\K[y]$ that are
contained in $A$. Let $u$ be the generator of $\mathfrak{a}$. It
is clear that $u$ has properties (\ref{u1}) and (\ref{u2}).
Suppose $f \in A \otimes B$ and $(\K[y] \otimes B)f \subseteq A
\otimes B$. Write $f = \sum_i a_i \otimes b_i$ for some $a_i \in
A$, $b_i \in B$, where the set $\{b_i\}$ is linearly independent
over $\K$. Let $q \in \K[y]$. Since $qf = \sum_i (q a_i) \otimes
b_i \in A \otimes B$, we have $q a_i \in A$ for each $i$. Element
$q$ was arbitrary, and so property (\ref{u2}) implies that $u$
divides each $a_i$ in $\K[y]$, and thus $u$ divides $f$ in $\K[y]
\otimes B$.
\end{proof}

An exponential map $\f:A \otimes B \to (A \otimes B)[t]$ is
uniquely extended to a homomorphism $\f: \Frac(A) \otimes B \to
\Frac((A \otimes B)[t])$ by setting $\f(a^{-1}) = \f(a)^{-1}$ for
all $a \in A \setminus 0$. Remark that this extension of $\f$
retains the property that $\varepsilon_0 \f$ is the identity map,
where $\varepsilon_0$ is evaluation at $t=0$.

\begin{lemma}\label{L:technical3}
Let $\f: A \otimes B \to (A \otimes B)[t]$ be an exponential map.
Let $\f$ also denote the unique extension $\Frac(A) \otimes B \to
\Frac((A \otimes B)[t])$. Then $\f(\K[y] \otimes B) \subseteq
(\K[y] \otimes B)[t]$.
\end{lemma}

\begin{proof}
It suffices to show that $\f(y) \in (\K[y] \otimes B)[t]$. Let $D$
be the locally finite iterative higher derivation associated to
$\f$. As mentioned earlier, $D$ has a unique extension to a higher
derivation on any given localization of $A \otimes B$. In
particular, there is a unique way to define the derivation on $y$.
Write $y = g h^{-1}$ for some $g, h \in A$. In any extension of
$D$, each $D^j(y)$ is found as some expression of elements from
$\{D^i(g)\}$ and $\{D^i(h)\}$, divided by some power of $h$.
Therefore, each $D^j(y)$ belongs to $\K(y) \otimes B$.

Since $y$ is integral over $A$, $\f(y)$ must be integral over $(A
\otimes B)[t]$. By viewing $(A \otimes B)[t]$ as a subring of
$B[y,t]$, we see that $\f(y)$ must belong to $\Frac(B)[y,t]$, i.e.
$\f(y) \in (\K[y] \otimes \Frac(B))[t]$. Thus we can restrict the
codomain of our extension of $\f$, so that we have $\f: \Frac(A)
\otimes B \to (\K[y] \otimes \Frac(B))[t]$. This extension of $\f$
defines an extension of $D$ with each $D^j(y)$ belonging to $\K[y]
\otimes \Frac(B)$. Combining this with our previous observation,
we have $D^j(y) \in \K[y] \otimes B$ for each $j$. Hence $\f(y)
\in (\K[y] \otimes B)[t]$.
\end{proof}

\begin{lemma}\label{L:technical4}
Let $u \in A$ be as in the statement of Lemma~\ref{L:technical2}.
Suppose $\f \in \EXP(A \otimes B)$ and let $D = \{D^i\}$ be the
locally finite iterative higher derivation associated to $\f$.
Then $u$ divides $D^n(u)$ in $\K[y] \otimes B$ for each natural
number $n$.
\end{lemma}

\begin{proof}
We will use induction to show that $(\K[y] \otimes B)D^n(u)
\subseteq A \otimes B$ and then appeal to
Lemma~\ref{L:technical2}. If $n=0$ then $D^n(u)=u$ and the result
is true by our choice of $u$. Suppose the result is true for
$D^l(u)$, where $0 \leq l < n$. Let $q \in \K[y] \otimes B$. By
Lemma~\ref{L:technical3}, we can extend $\f$ to a homomorphism
$\K[y] \otimes B \to (\K[y] \otimes B)[t]$. Thus we can write
$\f(q) = \sum_i D^i(q) t^i$, where $D^i(q) \in \K[y] \otimes B$
for all $i \geq 0$ and $D^0(q)=q$. Now $u, uq \in A \otimes B$,
and so $D^i(u), D^i(uq) \in A \otimes B$ for all $i \geq 0$. The
Leibniz rule yields
\[
D^n(uq) = \sum_{i=0}^{n-1} D^i(u) D^{n-i}(q) + D^n(u) q.
\]
By the induction hypothesis, each term in the summation belongs to
$A \otimes B$. Thus $D^n(u) q \in A \otimes B$ as well. We have
now verified that $(\K[y] \otimes B)D^n(u) \subseteq A \otimes B$
for each natural number $n$. The claim follows from
Lemma~\ref{L:technical2}.
\end{proof}

We are now in position to complete the proof. Since $A \nsubseteq
\ak(A \otimes B)$, there exists $\f \in \EXP(A \otimes B)$ for
which $A \cap (A \otimes B)^{\f} = \K$. Let $u \in A$ be as in the
statement of Lemma~\ref{L:technical2}. Let $n = \deg_{\f}(u)$.
Then $D^n(u) \in (A \otimes B)^{\f} \setminus 0$ by part
(\ref{factsc}) of Lemma~\ref{L:facts}. By
Lemma~\ref{L:technical4}, we have $D^n(u) = ug$ for some nonzero
$g \in \K[y] \otimes B$. Both $u$ and $ug^2$ are nonzero elements
of $A \otimes B$, and $u \cdot ug^2 = D^n(u)^2 \in (A \otimes
B)^{\f} \setminus 0$. By part (\ref{factsa}) of
Lemma~\ref{L:facts}, $u \in (A \otimes B)^{\f}$. Therefore, $u \in
\K$. Since $\K[y]u \subseteq A$, we have $A = \K[y]$. This
contradicts our assumption on $A$.


\end{document}